\documentclass[aps,11pt]{revtex4-1}

\begin{document}

\tighten

\draft

\title{Differential Operator Method of Finding A Particular Solution to
An Ordinary Nonhomogeneous Linear Differential Equation with
Constant Coefficients }

\author{ Wenfeng Chen\renewcommand{\thefootnote}{\dagger}\thanks{E-mail:
chenw3@sunyit.edu} }

\address{ Department of Mathematics and Physics, College of Arts and Sciences,
SUNY Polytechnic Institute, Utica, NY 13502, USA}


\begin{abstract}

\noindent We systematically introduce the idea of applying
differential operator method to find a particular solution of an
ordinary nonhomogeneous linear differential equation with constant
coefficients when the nonhomogeneous term is a polynomial
function, exponential function, sine function, cosine function or
any possible product of these functions. In particular, different
from the differential operator method introduced in literature, we
propose and highlight utilizing the definition of the inverse of
differential operator to determine a particular solution. We
suggest that this method should be introduced in textbooks and
widely used for determining a particular solution of an ordinary
nonhomogeneous linear differential equation with constant
coefficients in parallel to the method of undetermined
coefficients.


\end{abstract}

\maketitle

\vspace{3ex}

\section{Introduction}

\noindent  When we determine a particular solution of an ordinary
  nonhomogeneous linear differential equation of constant
coefficients with nonhomogeneous terms being a polynomial
function, an exponential function, a sine function, a cosine
function or any possible products of these functions, usually only
the method of undetermined coefficients is introduced in
undergraduate textbooks \cite{zill, boyce,snider}. The basic idea
of this method is first assuming the general form of a particular
solution with the coefficients undetermined, based on the
nonhomogeneous term, and then substituting the assumed solution to
determine the coefficients. The calculation process of this method
is very lengthy and time-consuming. In particular, in the case
that the assumed particular solution duplicates a solution to the
associated homogeneous equation, one must multiply a certain power
function to erase the duplication, and students usually get
confused and frustrated in applying this method to find a
particular solution of the equation.

Nevertheless, differential operator method provide a convenient
and effective method of finding a particular solution of an
ordinary nonhomogeneous linear differential equation of constant
coefficients with the nonhomogeneous terms being a polynomial
function, an exponential function, a sine function, a cosine
function or any possible products of these functions. The
efficiency of this method in determining a particular solution is
based on the following facts: exponential function is an
eigenfunction of $D\equiv d/dx$; sine and cosine functions are
eigenfunctions of second-order differential operator $D^2$; There
exists an exponential shift theorem when a polynomial of
differential operator acts on a product of exponential function
with another continuous function. Especially, if we apply the
inverse of differential operator acting on a differentiable
function, which is an integration over the function, this method
can greatly reduce the complication and arduousness  of the
calculation.

However, the introduction to this method only scatters partially
in some online lecture notes of differential equations
\cite{onlinenotes}. Hence, here we give a systematic introduction
to the method including both rigorous mathematical principles and
detailed  calculation techniques. In Sect.II, we introduce the
definitions of differential operator and its higher order versions
as well their actions on differentiable functions. We emphasize
the fact that exponential functions are eigenfunctions of
differential operator and both sine and cosine functions  are
eigenfunctions of differential operators of even orders. Further,
we define  the polynomial of differential operator and list its
properties when it acts on differentiable functions. We highlight
two remarkable properties including eigenvalue substitution rule
and exponential shift rule, which will play important roles in
solving differential equations. In addition, we introduce the
kernel of the polynomial operator and point out it is the solution
to the associated linear homogeneous differential equations. In
Sect. III, based on the Fundamental Theorem of Calculus, we define
the inverses of differential operator and its higher order
versions. Then according to the definition of a function of
differential operator given by the Taylor series we define the
inverse of a polynomial of differential operator and show some of
its properties when it acts  on differentiable functions. Sect. IV
contains the main content of the article. We use a large number of
examples to show how the differential operator method works
efficiently in determining a particular solution for a
nonhomogeneous linear differential equation with constant
coefficients when the nonhomogeneous term is a polynomial
function, an exponential function, a sine, a cosine function or
any possible product of these functions. Sect V is a a brief
summary on the application of differential operator method in
solving nonhomogeneous linear ordinary differential equations with
constant coefficients.

\section{Differential Operator and Action on A Differentiable Function }
\label{sect2}

\subsection{Definition of differential operator}

\noindent A differential operator $D$ acting a differentiable
function $y=f(x)$ on $R$ takes the form
\begin{eqnarray}
D=\frac{d}{dx} \label{eqdop}
\end{eqnarray}
The action of $D$ and its higher order versions $D^n$ on a at
least $n$-times differentiable function $y=f(x)$ is just to take
the derivatives of the function:
\begin{eqnarray}
Df(x)&=&\frac{df(x)}{dx}, \\ \nonumber
 D^2f(x)&=&
D[Df(x)]=\frac{d}{dx}\left[\frac{df(x)}{dx}\right]
=\frac{d^2f(x)}{dx^2}, \\ \nonumber &\vdots &
\nonumber\\
D^nf(x)&=&
D[D^{n-1}f(x)]=\frac{d}{dx}\left[\frac{d^{n-1}f(x)}{dx^{n-1}}\right]
=\frac{d^nf(x)}{dx^n} \label{opact}
\end{eqnarray}
For example, let
\begin{eqnarray}
y=x^3+x+3 e^{2x}-5\sin (3x)
\end{eqnarray}
then
\begin{eqnarray}
Dy=\frac{d}{dx}\left[x^3+x+3 e^{2x}-5\sin (3x)\right]
=3x^2+1+6e^{2x}+15 \cos (3x)
\end{eqnarray}

\subsection{Action of differential operator on elementary
functions and eigenvalue of differential operator}

The actions of differential operator on elementary functions
including exponential function, sine and cosine function, and
polynomial functions are listed as follows:

\begin{itemize}
\item[1.] \textbf{Exponential function}
\begin{eqnarray}
De^{\lambda x}=\frac{d}{dx}e^{\lambda x}=\lambda e^{\lambda x}
\end{eqnarray}
It shows $e^{\lambda x}$ is an eigenfunction of $D$ with
eigenvalue $\lambda$ in the representation space composed of
differentiable function on $R$.

Straightforwardly, for the higher order version $D^n$, when $n\geq
2$ is a positive integer, there exists
\begin{eqnarray}
D^ne^{\lambda x}=\frac{d^n}{dx^n}e^{\lambda
x}=\frac{d}{dx}\left[\frac{d}{dx}\left(\cdots\frac{d}{dx}e^{\lambda
x}\right)\right]=\lambda^n e^{\lambda x}
\end{eqnarray}

\item[2.] \textbf{Action of $D^2$ on sine and cosine functions}

\begin{eqnarray}
D^2 \sin (\beta x) = -{\beta}^2 \sin (\beta x), \nonumber\\
D^2 \cos (\beta x) = -{\beta}^2 \cos (\beta x) \label{actsc}
\end{eqnarray}
This shows both $\sin (\beta x)$ and $\cos (\beta x)$ are
eigenfunctions of $D^2$ with eigenvalue $-{\beta}^2$. Further,
consider the higher order version $D^{2n}=\left(D^2\right)^n$, we
have
\begin{eqnarray}
D^{2n}\sin (\beta x) =\left( D^2 \right)^n\sin (\beta x)
=D^2\left[D^2\left(\cdots D^2\sin \beta x\right)\right]
= \left(-{\beta}^2\right)^n \sin (\beta x),
 \nonumber\\
D^{2n}\cos (\beta x) =\left( D^2 \right)^n\cos (\beta x)
=D^2\left[D^2\left(\cdots D^2\cos \beta x\right)\right] =
\left(-{\beta}^2\right)^n \cos (\beta x)
\end{eqnarray}

\item[3.] \textbf{Action of $D$ on power function}

\begin{eqnarray}
Dx^k &=& \frac{d}{dx}x^k=kx^{k-1}, \nonumber\\
D^nx^k &=&\frac{d^n}{dx^n}x^k=k (k-1)\cdots (k-n+1)x^{k-n}
\end{eqnarray}
where $k\geq 1$ is a positive integer. Especially,
\begin{eqnarray}
D^nx^k =0~~\mbox{for $k<n$} \label{hactl}
\end{eqnarray}
For example,
\begin{eqnarray}
D^5x^3 =\frac{d^5}{dx^5}x^3=0
\end{eqnarray}

\end{itemize}

\subsection{Polynomial of Differential Operator $D$ and Property}

\noindent \textbf{Definition 1}: Let
\begin{eqnarray}
P_n(x)=a_n x^n+ a_{n-1}x^{n-1}+\cdots+ a_1 x +a_0
\end{eqnarray}
be a polynomial of real variable $x$ of degree $n$, where $a_0$,
$a_1$, $\cdots$, $a_n$ are real constants. Then
\begin{eqnarray}
P_n(D)=a_n D^n+ a_{n-1} D^{n-1}+a_1 D+\cdots +a_0
\label{polydoper}
\end{eqnarray}
 is called a polynomial of differential operator $D$ of degree $n$. In the
 representation space composed of at least $n$ times differentiable
 functions $f(x)$ on $R$,
 \begin{eqnarray}
P_n(D)=P_n\left(\frac{d}{dx}\right) =a_n \frac{d^n}{dx^n}+ a_{n-1}
\frac{d^{n-1}}{dx^{n-1}}+\cdots +a_1 D+a_0
\end{eqnarray}

The properties of a polynomial of differential operator $D$ are
listed as follows:

\begin{itemize}
\item[1.] \textbf{Linearity}

\textbf{Theorem 1:} Let $P_n(D)$ be a polynomial of differential
operator $D$ of degree $n$ defined in (\ref{polydoper}). Then
there exists
 \begin{eqnarray}
P_n(D)\left[ a f(x)+b g(x)\right]=a P_n(D)f(x)+b P_n(D)g(x)
\end{eqnarray}
where $a$ and $b$ are constants, and both $f(x)$ and $g(x)$ are
differentiable functions of at least $n$-times.

\item[2.] \textbf{Sum rule and product rule}

\noindent \textbf{Theorem 2:} Let $P_n(D)$ and $Q_m(D)$ be
polynomials of differential operator $D$ of degree $n$ and $m$,
respectively. Then there exist\\
 \textbf{(1). Sum rule}
 \begin{eqnarray}
\left[P_n(D)+Q_m (D)\right]f(x)=\left[Q_m (D)+P_n(D)\right]f(x) =
P_n(D)f(x)+Q_m(D)f(x)
\end{eqnarray}
where  $f(x)$ is a differentiable function of at least
$\mbox{max}(n,m)$ times.\\

\textbf{(2). Product  rule}
 \begin{eqnarray}
\left[P_n(D)Q_m (D)\right]g(x) &=& \left[Q_m (D)P_n(D)\right]g(x)
= P_n(D)\left[Q_m(D)g(x)\right]
\nonumber\\
&=& Q_m(D)\left[P_n(D)g(x)\right] \label{prodrule}
\end{eqnarray}
where  $f(x)$ is a differentiable function of at least
$\mbox{max}(n,m)$ times.

\textbf{Theorems 1} and \textbf{2} can be proved straightforwardly
with the explicit representation  forms of $P_n(D)$ and $Q_m(D)$
on the space of differentiable functions.

\item[3.] \textbf{Eigenvalue substitution rule}

\textbf{Theorem 3:} Let $P_n(D)$ be a polynomial of differential
operator $D$ of degree $n$. Then there exist
 \begin{eqnarray}
&(1).& ~~~ P_n(D) e^{\lambda x}=P_n(\lambda)e^{\lambda
x}=e^{\lambda x}\,P_n(\lambda) \label{eigensub1}
\\
&(2).& ~~~ P_n(D^2) \sin\beta x =P_n(-\beta^2)\sin\beta x
=\sin\beta
x\,P_n(-\beta^2), \nonumber\\
& & ~~~P_n(D^2) \cos\beta x =P_n(-\beta^2)\cos\beta x =\cos \beta
x\, P_n(-\beta^2) \label{eigensub2}
\end{eqnarray}
where $\lambda$ and $\beta$ are real constants.

\textbf{Proof:} (1). First we have
 \begin{eqnarray}
P_n(D)e^{\lambda x}&=&P_n\left(\frac{d}{dx}\right)e^{\lambda x}
=\left(a_n \frac{d^n}{dx^n}+ a_{n-1}
\frac{d^{n-1}}{dx^{n-1}}+\cdots +a_1
D+a_0\right)e^{\lambda x}\nonumber\\
&=& \left(a_n \lambda^n+ a_{n-1} \lambda^{n-1}+\cdots +a_1 \lambda
+a_0\right)e^{\lambda x}=P_n(\lambda)e^{\lambda x}=e^{\lambda
x}P_n(\lambda)
\end{eqnarray}

(2). First applying Eq.(\ref{actsc}), we obtain
\begin{eqnarray}
P_n(D^2) \sin\beta x &=& \left[a_n (D^2)^n+a_{n-1}
(D^2)^{n-1}+\cdots+a_1 D^2+a_0 \right]\sin\beta x \nonumber\\
&=&\left[a_n (-\beta^2)^n+a_{n-1} (-\beta^2)^{n-1}+\cdots+a_1
(-\beta^2)+a_0 \right]\sin\beta x\nonumber\\
&=& P_n(-\beta^2) \sin\beta x =\sin\beta x\,P_n(-\beta^2)
\end{eqnarray}
Similarly there is the result for $P_n(D^2) \cos\beta x$.

 \item[4.] \textbf{Exponential shift rule}

\textbf{Theorem 4:} Let $P_n(D)$ be a polynomial of differential
operator $D$ of degree $n$. Then there exists
 \begin{eqnarray}
 D^n\left[ e^{\lambda x}f(x)\right]&=& e^{\lambda x} \left(D+\lambda
 \right)^n f(x), \label{dexpshift} \\
P_n(D)\left[ e^{\lambda x} f(x)\right] &=& e^{\lambda x} P_n
\left(D+\lambda \right)f(x) \label{polyshift}
\end{eqnarray}
where $f(x)$ is a differentiable function of at least $n$ times on
$R$.

\textbf{Proof:} \textbf{(1).} We use the method of mathematical
induction to prove the result (\ref{dexpshift}).

First, the result arises for $n=1$:
\begin{eqnarray}
 D\left[ e^{\lambda x}f(x)\right]&=&  \frac{d}{dx}\left[ e^{\lambda x}f(x)\right]
 =e^{\lambda x}\left(\lambda +\frac{d}{dx}\right) f(x)=e^{\lambda
 x}\left( D+\lambda\right) f(x)
\end{eqnarray}
Second, assume the result holds  for $n=k$:
\begin{eqnarray}
D^k\left[ e^{\lambda x}f(x)\right]= e^{\lambda x} \left(D+\lambda
 \right)^k f(x)
\end{eqnarray}
Then
\begin{eqnarray}
D^{k+1}\left[ e^{\lambda x}f(x)\right]&=&
D\left[D^k\left(e^{\lambda x} f(x)\right)\right]=D\left[e^{\lambda
x} \left(D+\lambda
 \right)^k f(x)\right]\nonumber\\
 &=& e^{\lambda x}(D+\lambda )\left[ \left(D+\lambda
 \right)^k f(x)\right]=e^{\lambda x} \left(D+\lambda
 \right)^{k+1}f(x)
\end{eqnarray}
Therefore, we have the general  result (\ref{dexpshift}).

\textbf{(2).} The result (\ref{polyshift}) follows
straightforwardly from (\ref{dexpshift}):
\begin{eqnarray}
&& P_n(D)\left[ e^{\lambda x} f(x)\right]= \left[a_n D^n+a_{n-1}
D^{n-1}+\cdots+a_1 D+a_0 \right]\left[ e^{\lambda x}
f(x)\right]\nonumber\\
& = &e^{\lambda x}\left[ a_n (D+\lambda)^n+a_{n-1}
(D+\lambda)^{n-1}+\cdots+a_1 (D+\lambda)+a_0 \right]f(x)\nonumber\\
&=& e^{\lambda x} P_n (D+\lambda)f(x)
\end{eqnarray}

\end{itemize}

\subsection{Kernel of Differential Operator}

\noindent According to the general definition on the kernel of an
operator, the kernel of differential operator $D^k$, $k=1,2,
\cdots, $ in the space of differentiable real functions is defined
as follows:

\noindent \textbf{Definition 3}:
\begin{eqnarray}
\mbox{ker}D^k=\left\{\left.f(x)\right| D^k f(x)=0  \right\}
\end{eqnarray}
Due to Eq.(\ref{hactl}), $\mbox{ker}D^k$ is the set of polynomial
functions of degree at most $k-1$,
\begin{eqnarray}
\mbox{ker}D^k &=&\mbox{span}\left\{1, x , \cdots,
x^{k-1}\right\}\nonumber\\
 &=&\left\{\left.c_{k-1}x^{k-1}+c_{k-2}x^{k-2}+\cdots+c_1x+c_0\right|
c_0, c_1, \cdots, c_k~\mbox{are constants} \right\}
\end{eqnarray}

Further, we can write down the kernel of the  polynomial $P_n(D)$
of differential operator $D$, which is the solution space of a
linear $n$-th order ordinary homogeneous differential equation
with constant coefficients,
\begin{eqnarray}
P_n(D)y(x)=a_ny^{(n)}(x)+a_{n-1}y^{(n-1)}(x)+\cdots+ a_1y^\prime
(x)+a_0 y(x)=0
\end{eqnarray}

The kernel $\mbox{ker}P_n(D)$ takes different forms depending on
the form of $P_n(D)$. Applying the exponential shift rule
(\ref{dexpshift}) and eigenvalue substitution rule, we can easily
find the following results:

\noindent \textbf{Theorem 5}:
\begin{itemize}
\item[(1).] If $P_n(D)$ is product of $n$ distinct linear factors
(for simplicity, we take $a_n=1$ for now and later),
\begin{eqnarray}
P_n(D)&=&(D-r_1) (D-r_2)\cdots
(D-r_n)=\prod_{i=1}^n(D-r_i),\nonumber\\
&& r_1\neq r_2\neq \cdots \neq r_n
\end{eqnarray}
then
\begin{eqnarray}
\mbox{ker}P_n(D)&=& \mbox{span}\left\{e^{r_1x}, e^{r_2x}, \cdots,
e^{r_nx} \right\}\nonumber \\
&=&\left\{\left.\sum_{i=1}^n c_i
e^{r_ix}\right| c_1, c_2, \cdots, c_n~\mbox{are constants}\right\}
\end{eqnarray}

\item[(2).] If $P_n(D)$ is product of $k$ distinct repeated linear
factors,
\begin{eqnarray}
P_n(D)&=&(D-r_1)^{n_1} (D-r_2)^{n_2}\cdots
(D-r_k)^{n_k}=\prod_{i=1}^k (D-r_i)^{n_i},\nonumber\\
&& r_1\neq r_2\neq \cdots \neq r_k \geq 1, ~~\sum_{i=1}^k n_i=n
\end{eqnarray}
then
\begin{eqnarray}
\mbox{ker}P_n(D)= \mbox{span}\left\{xe^{r_1x}, x^2e^{r_1x},\cdots,
x^{n_1-1}e^{r_1x}; \cdots; xe^{r_k x}, x^2e^{r_k x}, \cdots,
x^{n_k-1}e^{r_k x}\right\}
\end{eqnarray}

\item[(3).] If $P_n(D)$ is product of $k$ distinct quadratic
factors, and in this case there must be $n=2k$,
\begin{eqnarray}
P_n(D)&=&\left[(D-\alpha_1)^2+\beta_1^2\right]
\left[(D-\alpha_2)^2+\beta_2^2\right] \cdots
\left[(D-\alpha_k)^2+\beta_k^2\right]\nonumber\\
 &=& \prod_{i=1}^k \left[(D-\alpha_i)^2+\beta_i^2\right],~~
  \alpha_i\neq \alpha_j, ~\mbox{or}~ \beta_i\neq \beta_j~\mbox{for}~i\neq j
\end{eqnarray}
then
\begin{eqnarray}
\mbox{ker}P_n(D)= \mbox{span}\left\{e^{\alpha_1x}\cos(\beta_1 x),
e^{\alpha_1x}\sin(\beta_1 x); \cdots; e^{\alpha_kx}\cos(\beta_k
x), e^{\alpha_kx}\sin(\beta_k x) \right\}
\end{eqnarray}

\item[(4).] If $P_n(D)$ is product of $p$ repeated distinct
quadratic factors,
\begin{eqnarray}
P_n(D)&=&\left[(D-\alpha_1)^2+\beta_1^2\right]^{n_1}
\left[(D-\alpha_2)^2+\beta_2^2\right]^{n_2} \cdots
\left[(D-\alpha_p)^2+\beta_p^2\right]^{n_p}\nonumber\\
 &=& \prod_{i=1}^p \left[(D-\alpha_i)^2+\beta_i^2\right]^{n_i},~~
  \alpha_i\neq \alpha_j, ~\mbox{or}~ \beta_i\neq \beta_j~\mbox{for}~i\neq j
\end{eqnarray}
then
\begin{eqnarray}
\mbox{ker}P_n(D)&=& \mbox{span}\left\{e^{\alpha_1x}\cos(\beta_1
x), e^{\alpha_1x}\sin(\beta_1 x), x e^{\alpha_1x}\cos(\beta_1 x),
x e^{\alpha_1x}\sin(\beta_1 x), \right. \nonumber\\
&& \cdots, x^{n_1-1} e^{\alpha_1x}\cos(\beta_1 x),
x^{n_1-1} e^{\alpha_1x}\sin(\beta_1 x); \nonumber  \\
&& e^{\alpha_2x}\cos(\beta_2 x), e^{\alpha_2 x}\sin(\beta_2 x), x
e^{\alpha_2 x}\cos(\beta_2 x), x e^{\alpha_2 x}\sin(\beta_2 x), \nonumber \\
&& \cdots, x^{n_2-1} e^{\alpha_2 x}\cos(\beta_2 x), x^{n_2-1}
e^{\alpha_2 x}\sin(\beta_2 x);\nonumber \\
&&  \cdots; e^{\alpha_p x}\cos(\beta_p x), e^{\alpha_p
x}\sin(\beta_p x), x e^{\alpha_p x}\cos(\beta_p x), x e^{\alpha_p
x}\sin(\beta_p x), \nonumber \\
&& \left. \cdots, x^{n_p-1} e^{\alpha_p x}\cos(\beta_p x),
x^{n_p-1} e^{\alpha_p x}\sin(\beta_p x) \right\}
\end{eqnarray}

\end{itemize}

\noindent If $P_n(D)$ is a mixture of the above four cases, so are
the corresponding terms from each of the kernels. For example, let
\begin{eqnarray}
P_{12}(D)= \left(D-2\right)\left(D-5\right)^3
\left[\left(D+3\right)^2+4 \right]
\left[\left(D-7\right)^2+16\right]^4
\end{eqnarray}
then
\begin{eqnarray}
\mbox{ker}P_{12}(D)&=& \mbox{span}\left\{e^{2x}; e^{5x}, x e^{5x},
x^2 e^{5x}; e^{-3x}\cos(2x), e^{-3x}\sin (2x); e^{7x} \cos (4x),
e^{7x} \sin (4x), \right.
\nonumber\\
&& \left. x e^{7x} \cos (4x), xe^{7x} \sin (4x), x^2e^{7x} \cos
(4x), x^2 e^{7x} \sin (4x), x^3e^{7x} \cos (4x), x^3 e^{7x} \sin
(4x) \right\}\nonumber\\
&=& a_0e^{2x}+\left(b_0+b_1x+b_2x^2\right)e^{5x}+\left(c_1\cos
2x+c_2\sin 2x\right)e^{-3x}\nonumber\\
&+&\left(d_0+d_1x+d_2 x^2+d_3 x^3 \right)e^{7x}\cos
4x+\left(e_0+e_1x+e_2 x^2+e_3 x^3 \right)e^{7x}\sin 4x
\end{eqnarray}
where $a_0$, $b_0$, $b_1$, $b_2$, $c_1$, $c_2$, $d_0$, $\cdots$,
$d_3$ and $e_0$, $\cdots$, $e_3$  are real numbers.

\section{Inverse of Differentia Operator and Action on A Continuous Function}

\subsection{Definition of the Inverse of Differential Operator}

\noindent The inverse of differential operator $D$ can be defined
according to the Fundamental Theorem of Calculus. Because there
exists
\begin{eqnarray}
\frac{d}{dx}\int_{x_0}^x f(x_1)dx_1=D\int_{x_0}^x f(x_1)dx_1=f(x)
\end{eqnarray}
where $f(x)$ is a continuous function on a finite interval $[a,b]$
and $x_0$ is an arbitrary constant on $[a, b]$, therefore, the
action of the inverse $D^{-1}$ of the differential operator $D$ on
a continuous function $f(x)$ can be defined in terms of the
following definite integral:
\begin{eqnarray}
D^{-1}f(x)=\int_{x_0}^x f(x_1)dx_1 \label{inverd1}
\end{eqnarray}
Further, we can successively define the action of higher order
inverse differential operator on a continuous function on a finite
interval:
\begin{eqnarray}
D^{-2}f(x)&=&\left(D^{-1}\right)^2f(x)=D^{-1}\left[
D^{-1}f(x)\right]=
\int_{x_0}^x dx_1\int_{x_0}^{x_1}f(x_2)dx_2, \label{inverd2}\\[2mm]
D^{-3}f(x)&=&\left(D^{-1}\right)^3f(x) = \int_{x_0}^x
dx_1\int_{x_0}^{x_1}dx_2\int_{x_0}^{x_2} f(x_3)dx_3,
\label{inverd3}\\[2mm]
&\vdots & \nonumber \\
D^{-n}f(x)&=&\left(D^{-1}\right)^nf(x) = \int_{x_0}^x
dx_1\int_{x_0}^{x_1}dx_2\cdots \int_{x_0}^{x_{n-1}} f(x_n)dx_n
\label{inverdn}
\end{eqnarray}

The justification of these definitions can be demonstrated by the
following examples:
\begin{eqnarray}
&1.& ~~ D^{-1}\cos (3x)=\int_{x_0}^x \cos (3
x_1)dx_1=\frac{1}{3}\left[ \sin (3x)-\sin
(3x_0)\right]\nonumber\\
 &&=\frac{1}{3}\sin (3x)+c= \frac{1}{3}\sin (3x)+\mbox{ker}D\\
&2.& ~~ D^{-2} x =\int_{x_0}^x \int_{x_0}^{x_1}dx_2 x_2
dx_2=\int_{x_0}^x \frac{1}{2}\left( x_1^2- x_0^2
\right)=\frac{1}{6}x^3-\frac{1}{2}x_0^2 x+\frac{1}{3}x_0^3\nonumber\\
&&~~ =\frac{1}{6}x^3+c_1x+c_2=\frac{1}{6}x^3+\mbox{ker}D^2
\end{eqnarray}

\subsection{Function of Differential Operator}

\noindent As a generalization of the polynomial $P_n(D)$, a
function of the differential operator $D$ can be defined in a way
similar to the definition of a square matrix \cite{dpool}.

\noindent \textbf{Definition 2}: Let $f(x)$ be a uniformly
convergent function on $R$. Then a function of $f(D)$ of the
differential $D$ is defined in terms of the Taylor series
expansion of $f(x)$ about $x=0$:
\begin{eqnarray}
f(D)=\sum_{n=0}^{\infty}\frac{1}{n!}f^{(n)}(0)D^n
=\sum_{n=0}^{\infty}\left.\frac{1}{n!}\frac{d^nf(x)}{dx^n}\right|_{x=0}
D^n=\sum_{n=0}^{\infty}\left.\frac{1}{n!}\frac{d^nf(x)}{dx^n}\right|_{x=0}
\frac{d^n}{dx^n} \label{taylordef}
\end{eqnarray}
For examples, two typical functions of $D$ are listed as follows:
\begin{eqnarray}
&(1).& ~~\frac{1}{1-D}=\sum_{n=0}^n D^n=1+D+D^2+\cdots+
D^n+\cdots \\
 &(2).& ~~ e^D=\sum_{n=0}^{\infty}\frac{1}{n!} D^n
 =1+D+\frac{1}{2!}D^2+\cdots+\frac{1}{n!} D^n+\cdots
\end{eqnarray}

\subsection{Inverse of Polynomial of Differential Operator}

\noindent Using the Taylor series representation (\ref{taylordef})
on a function of differential operator $D$, the product rule
(\ref{prodrule}) of polynomials of differential operator $D$, and
the definitions (\ref{inverd1}), (\ref{inverd2}), (\ref{inverd3}),
(\ref{inverdn}) of the action of the inverse operator $D^{-k}$ on
a continuous function, we can define the inverse of a polynomial
operator $P_n(D)$ as follows:

\noindent \textbf{Definition 4}: Let
\begin{eqnarray}
P_n(D)f(x)=\left(a_nD^n+a_{n-1}D^{n-1}+\cdots+a_1D+a_0\right)f(x)=g(x)
\end{eqnarray}
where both $f(x)$ and $g(x)$ are differentiable functions.
\begin{itemize}
\item[(1).] If $a_0\neq 0$, then
\begin{eqnarray}
f(x)&=& \left[P_n(D) \right]^{-1}g(x)= \frac{1}{P_n(D)}g(x) =
\frac{1}{a_nD^n+a_{n-1}D^{n-1}+\cdots+a_1D+a_0}\,g(x) \nonumber\\
&=&
\sum_{p=0}^{\infty}\left(-1\right)^p\frac{1}{a_0}\left(\frac{a_n}{a_0}D^n+\cdots
+\frac{a_1}{a_0}D \right)^p\, g(x) \label{def4part1}
\end{eqnarray}

\item[(2).] If all $a_0=a_1=\cdots a_{k-1}= 0$ ($1\leq k \leq n$)
and $a_k\neq 0$ , then
\begin{eqnarray}
f(x)&=& \left[P_n(D) \right]^{-1}g(x)= \frac{1}{P_n(D)}g(x) =
\frac{1}{a_nD^n+a_{n-1}D^{n-1}+\cdots+a_kD^k}\,g(x) \nonumber\\
&=&\frac{1}{a_nD^n+a_{n-1}D^{n-1}+\cdots+a_kD^k}\,g(x)\nonumber\\
&=& \frac{1}{a_nD^{n-k}+a_{n-1}D^{n-k-1}+\cdots+a_k}\,D^{-k} g(x)
\label{def4part2}
\end{eqnarray}

\end{itemize}

\noindent The definition of $\left[P_n(D) \right]^{-1}$ and the
product rule (\ref{prodrule}) yields the following property

\noindent \textbf{Corollary 1}: Let $P_n(D)$ and $Q_m(D)$ be two
polynomials of $D$ of degree $n$ and $m$, respectively. Then
\begin{eqnarray}
&& \frac{1}{P_n(D) Q_m(D)} f(x) =\frac{1}{Q_m(D)P_n(D) }
f(x)\nonumber\\
& =& \frac{1}{P_n(D)}\left[\frac{1}{ Q_m(D)} \right] f(x)
=\frac{1}{Q_m(D)}\left[\frac{1}{P_n(D) } \right] f(x)
\end{eqnarray}
where $f(x)$ is a differentiable function on a certain interval.

\noindent \textbf{Corollary 2}:  Let $P_n(D)$ and $Q_m(D)$ be two
polynomials of $D$ of degree $n$ and $m$, respectively. Then
\begin{eqnarray}
&& \frac{1}{P_n(D)} f(x) =\frac{1}{P_n(D)Q_m(D)} \left[
Q_m(D)f(x)\right]=Q_m(D) \left[\frac{1}{P_n(D)Q_m(D)} f(x)\right]
\end{eqnarray}
where $f(x)$ is a differentiable function on a certain interval.

\section{Differential Operator Method of Solving Nonhomogeneous Linear Ordinary
Differential Equation with Constant Coefficients}

 \noindent The general form of an $n$-th order nonhomogeneous linear ordinary
 differential equation with constant coefficients takes the following form:
\begin{eqnarray}
a_n\frac{d^ny}{dx^n}+
a_{n-1}\frac{d^{n-1}y}{dx^{n-1}}+\cdots+a_1\frac{dy}{dx}+a_0y=g(x)
\end{eqnarray}
where $a_0$, $a_1$, $\cdots$, $a_n$ are constants. Expressed in
terms of differential operator, the equation reads
\begin{eqnarray}
\left(a_n D^n+ a_{n-1}D^{n-1}+\cdots+a_1D+a_0\right)y=P_n(D)y=g(x)
\end{eqnarray}
Then formally a particular solution of the equation is
\begin{eqnarray}
Y(x)=\frac{1}{P_n(D)}g(x) \label{formsolu}
\end{eqnarray}

Due to the exponential shift rule (\ref{polyshift}) and the
eigenvalue substitution rule (\ref{eigensub1}), (\ref{eigensub2})
as well as the feature (\ref{hactl}), this approach, like the
method of undetermined coefficients, works only when $g(x)$ is a
polynomial function $P_k(x)$ of degree $k$,  an exponential
function $e^{\alpha x}$, $\sin (\beta x)$, $\cos (\beta x)$ or any
possible product of these functions. In the following we discuss
each case of $g(x)$

\subsection{The case $g(x)$ is an exponential function}

\noindent We first give the following results for $g(x)=A
e^{\alpha x}$ in order to show how the differential operator
approach works in this case, where $A$ is a constant.



\noindent \textbf{Theorem 6:} A particular solution to a
differential equation
\begin{eqnarray}
P_n(D) y(x)=A e^{\alpha x} \label{expdiffeq}
\end{eqnarray}
reads as follows:
\begin{itemize}
\item[1.] If $P_n(\alpha)\neq 0$, then
\begin{eqnarray}
 Y(x)=\frac{1}{P_n(D)}\left(A e^{\alpha x}\right)
 = \frac{1}{P_n(\alpha)}\left(A e^{\alpha x}\right)\label{expdeqsolu1}
\end{eqnarray}
\item[2.] If $P_n(\alpha)= 0$, i.e., $P_n(D)=(D-\alpha)^k
P_{n-k}(D)$, $1\leq k \leq n$ then
\begin{eqnarray}
 Y(x)=\frac{1}{P_n(D)}\left(A e^{\alpha x}\right)
 =\frac{A}{P_{n-k}(\alpha)}\left(\frac{1}{k!}x^k+\mbox{ker}D^k\right)e^{\alpha x}
  \label{expdeqsolu2}
\end{eqnarray}

\noindent\textbf{Proof:} According to the rule of eigenvalue
substitution (\ref{eigensub1}), we have
\begin{eqnarray}
 P_{n}(D)A e^{\alpha x}=A P_{n}(\alpha )e^{\alpha x} \label{eigensubeqrep}
\end{eqnarray}
\item[(1).] Since $P_n(\alpha)\neq 0$, we divide the equation
(\ref{eigensubeqrep}) by $P_n(\alpha)$, then
\begin{eqnarray}
\frac{1}{P_n(\alpha)}\left[P_n(D)\left( Ae^{\alpha
x}\right)\right]=P_n(D)\left(\frac{ Ae^{\alpha x}}{P_n
(\alpha)}\right)=A e^{\alpha x}
\end{eqnarray}
In comparison with the differential equation (\ref{expdiffeq}), it
shows $Y(x)=A e^{\alpha x}/P_n(\alpha)$.

\item[(2).] Since $P_n(\alpha)=0$, the polynomial $P_n(D)$ must
take the following form for a certain $1\leq k \leq n$,
\begin{eqnarray}
P_n(D)=P_{n-k}(D) \left(D-\alpha \right)^k, ~~P_{n-k}(\alpha)\neq
0,
\end{eqnarray}
and the equation becomes
\begin{eqnarray}
P_n(D)y(x)=(D-\alpha)^k P_{n-k}(D) y(x)=A e^{\alpha x}
\label{equform2}
\end{eqnarray}
Because
\begin{eqnarray}
&& (D-\alpha)^k P_{n-k}(D) \left[\left(\frac{1}{k!}x^k
+\mbox{ker}D^k\right)e^{\alpha x}\right]\nonumber\\
& =& P_{n-k}(\alpha)(D-\alpha)^k\left[
 e^{\alpha x}  \left(\frac{1}{k!}x^k
+\mbox{ker}D^k\right)\right] \nonumber\\
&=& P_{n-k}(\alpha) e^{\alpha x}D^k \left(\frac{1}{k!}x^k
+\mbox{ker}D^k\right) =P_{n-k}(\alpha) e^{\alpha x}
\end{eqnarray}
therefore, $y(x)$ given in (\ref{expdeqsolu2}) is the solution to
the differential equation.

\end{itemize}

\noindent \textbf{Theorem 6} shows a formal method of finding a
particular solution of differential equation when the
nonhomogeneous term is an exponential function,
\begin{eqnarray}
Y(x)&=& \frac{1}{P_n(D)}\left(A e^{\alpha
x}\right)=\frac{1}{P_{n-k}(D) \left(D-\alpha\right)^k)}\left(A
e^{\alpha x}\right)\nonumber\\
&=& \frac{A e^{\alpha x}}{P_{n-k}(\alpha)}D^{-k}\,1 =\frac{A x^k
e^{\alpha x}}{k!P_{n-k}(\alpha)}
\end{eqnarray}

The following two examples shows how this method works.
\begin{itemize}
\item[(1).] A particular solution to the differential equation
\begin{eqnarray}
3y^{\prime\prime}-2y^{\prime}+6y=5 e^{3x}
\end{eqnarray}
is
\begin{eqnarray}
Y(x)=\frac{1}{3 D^2-2D+8}
\left(5e^{3x}\right)=\frac{5e^{3x}}{3\cdot 3^2-2 \cdot
3+8}=\frac{5}{29}e^{3x}
\end{eqnarray}

\item[(2).] A particular solution to the differential equation
\begin{eqnarray}
 (D-1)(D+5)(D-2)^3 y(x)=3 e^{2x}
\end{eqnarray}
is
\begin{eqnarray}
Y(x)&=&\frac{1}{(D-1)(D+5)(D-2)^3}\left(3 e^{2x}\right)
 =\frac{3e^{2x}}{(2-1)(2+5)} D^{-3}1\nonumber\\
 &=& \frac{3e^{2x}}{7}
 \frac{1}{3!}x^3=\frac{1}{14}x^3e^{2x}
\end{eqnarray}

\end{itemize}


\noindent A straightforward result comes from \textbf{Theorem 6},
which is called the Exponential Input Theorem.

\noindent \textbf{Corollary 3 (Exponential Input Theorem)}:

Let $P_n(D)$ be a polynomial of $D$ of degree $n$. Then the
nonhomogeneous linear differential equation (\ref{expdiffeq}) has
a particular solution of the following form:
\begin{eqnarray}
Y(x)=\left\{\begin{array}{ll} \displaystyle \frac{Ae^{\alpha
x}}{P_n(\alpha)}, &
\mbox{if}~ P_n(\alpha)\neq 0 \\[2mm]
\displaystyle \frac{Ax e^{\alpha x}}{P_n^\prime (\alpha)}, &
\mbox{if}~ P_n(\alpha)=0~\mbox{but}~ P_n^\prime (\alpha) \neq
0\\[2mm]
\displaystyle \frac{Ax^2 e^{\alpha x}}{P_n^{\prime\prime}
(\alpha)}, &
\mbox{if}~ P_n(\alpha)=P_n^\prime (\alpha)=0
~\mbox{but}~P_n^{\prime\prime} (\alpha) \neq 0\\
 ~~ \vdots & \\[2mm]
\displaystyle \frac{Ax^k e^{\alpha x}}{P_n^{(k)} (\alpha)}, &
\mbox{if}~ P_n(\alpha)=P_n^\prime (\alpha)=\cdots =P_n^{(k-1)}
(\alpha)= 0 ~\mbox{but}~P_n^{(k)} (\alpha) \neq 0
\end{array} \right.
\end{eqnarray}

\vspace{2mm}

\noindent \textbf{Proof:} We still start from
Eq.(\ref{eigensubeqrep}).

\begin{itemize}
\item[(1).] For the case $P_n(\alpha)\neq 0$, the proof is the
same as that for \textbf{Theorem 6}.

\item[(2).] If $P_n(\alpha)= 0$, while $P_n^\prime (\alpha)\neq 0$
we first take the derivative on the
 the equation
(\ref{eigensubeqrep}) with respect to $\alpha$,  then
\begin{eqnarray}
&& \frac{d}{d\alpha} \left[P_n(D)\left( Ae^{\alpha
x}\right)\right]= \frac{d}{d\alpha}\left[ AP_n (\alpha)e^{\alpha
x}\right],\nonumber\\
&& P_n(D)\left( Ax e^{\alpha x}\right)=A\left[P_n^\prime (\alpha)
e^{\alpha x}+P_n(\alpha)\left( xe^{\alpha
x}\right)\right]=AP_n^\prime (\alpha)e^{\alpha x}
\end{eqnarray}
We divide the above equation by $P_n^\prime (\alpha)$ and obtain
\begin{eqnarray}
\frac{1}{P_n^\prime (\alpha)}\left[P_n(D)\left( Ax e^{\alpha
x}\right)\right]=P_n(D)\left[\frac{Ax e^{\alpha x}}{P_n^\prime
(\alpha)} \right]=Ae^{\alpha x}
\end{eqnarray}
Therefore, a particular solution of the  differential equation
(\ref{expdiffeq}) is $y(x)=A x e^{\alpha x}/P_n^\prime (\alpha)$.

\item[(3).] If $P_n(\alpha)=P_n^\prime (\alpha)= 0$, while
$P_n^{\prime\prime}(\alpha)\neq 0$, we take the derivative the
equation (\ref{eigensubeqrep}) twice with respect to $\alpha$,
\begin{eqnarray}
&& \frac{d^2}{d\alpha^2} \left[P_n(D)\left( Ae^{\alpha
x}\right)\right]= \frac{d^2}{d\alpha^2}\left[ AP_n
(\alpha)e^{\alpha
x}\right],\nonumber\\
&& P_n(D)\left( Ax^2 e^{\alpha x}\right)=A\left[P_n^{\prime\prime}
(\alpha) e^{\alpha x}+2P_n^\prime (\alpha) xe^{\alpha
x}+P_n^\prime (\alpha) x^2e^{\alpha x}\right]=AP_n^{\prime\prime}
(\alpha)e^{\alpha x}
\end{eqnarray}
Then dividing the above equation by $P_n^\prime (\alpha)$, we have
\begin{eqnarray}
P_n(D)\left[\frac{Ax2 e^{\alpha x}}{P_n^{\prime\prime} (\alpha)}
\right]=Ae^{\alpha x}
\end{eqnarray}
Hence in this case the solution to the differential equation is
$y(x)=A x^2 e^{\alpha x}/P_n^{\prime\prime} (\alpha)$
(\ref{expdiffeq}).

\item[(4).] For the general case,  $P_n(\alpha)=P_n^\prime
(\alpha)= \cdots= P^{(k-1)}_n (\alpha)=0$, while
$P_n^{(k)}(\alpha)\neq 0$, we take the derivative on the equation
(\ref{eigensubeqrep}) $k$ times with respect to $\alpha$,
\begin{eqnarray}
 \frac{d^k}{d\alpha^k} \left[P_n(D)\left( Ae^{\alpha
x}\right)\right]= \frac{d^k}{d\alpha^k}\left[ AP_n
(\alpha)e^{\alpha x}\right]
\end{eqnarray}
With an application of the product rule for higher order
derivative,
\begin{eqnarray}
\frac{d^k}{dx^k}\left[ f(x)g(x)\right]
=\sum_{p=0}^k\left(\begin{array}{c} k \\ p\end{array}\right)
f^{(k-p)}(x) g^{(p)}(x),
\end{eqnarray}
it yields
\begin{eqnarray}
&& P_n(D)\left( Ax^k e^{\alpha x}\right)=
A\sum_{p=0}^k\left(\begin{array}{c} k \\
p\end{array}\right)P_n^{(k-p)}\left(\alpha\right)x^p e^{\alpha
x}\nonumber\\
&=& A\left[ P_n^{(k)} (\alpha) e^{\alpha x}+k P_n^{(k-1)} (\alpha)
xe^{\alpha x}+\frac{k(k-1)}{2}P_n^{(k-2)} (\alpha) x^2e^{\alpha
x}+ \cdots +  P_n (\alpha) x^k e^{\alpha
x}\right]\nonumber\\
&=& AP_n^{(k)} (\alpha)e^{\alpha x}
\end{eqnarray}
Therefore, a particular solution of the equation in this case is
$y(x)= A x^k e^{\alpha x}/P_n^{(k)} (\alpha)$.

\end{itemize}

\noindent One can also apply \textbf{Corollary 3} to find a
particular solution of a differential equation. For example,
determine a particular solution of the equation
\begin{eqnarray}
(D-2) (D-4)^3y=5 e^{4t}
\end{eqnarray}
Since
\begin{eqnarray}
&& P^{\prime}(D)=(D-4)^2 (4D-10),
~P^{\prime\prime}(D)=(D-4)(12D-36), ~P^{(3)}(D)=12 (D-7),
\nonumber\\
&& P(4)=P^\prime (4)=P^{\prime\prime}(4)=0, ~~P^{(3)}(4)\neq 0
\end{eqnarray}
so a particular solution of the equation is
\begin{eqnarray}
Y(x)=\frac{5 x^3 e^{4x}}{P^{(3)}(4)}=-\frac{5}{36}x^3e^{4x}
\end{eqnarray}

\subsection{The case $g(x)$ is a polynomial function $P_k(x)$}

\noindent In this case, due to the fact (\ref{hactl}) it is most
convenient to use \textbf{Definition 4} ( see
Eqs.(\ref{def4part1}) and (\ref{def4part1})) to find a particular
solution of a differential equation.  We use the following two
examples to illustrate the idea.

\noindent \textbf{Example 1:} Determine a particular solution of
the equation
\begin{eqnarray}
y^{\prime\prime\prime}-5y^{\prime\prime}+3y^{\prime}+2y= x^2+3x-2
\end{eqnarray}
The differential operator form of the equation is
\begin{eqnarray}
\left(D^3-5D^2+3D+2\right)y= 2x^3+4x^2-6x+5
\end{eqnarray}
and consequently, a particular solution of the equation is
\begin{eqnarray}
Y(x)&=&\frac{1}{D^3-5D^2+3D+2}\left( 2x^3+4x^2-6x+5
\right)\nonumber\\
 &=&\frac{1}{2}\frac{1}{1+(D^3-5D^2+3D)/2}\left(
2x^3+4x^2-6x+5 \right) \nonumber\\
&=& \frac{1}{2}\left[
1-\frac{1}{2}(D^3-5D^2+3D)+\frac{1}{4}(-5D^2+3D)^2-\frac{1}{8}
(3D)^3\right]\left( 2x^3+4x^2-6x+5 \right)\nonumber\\
&=&
\frac{1}{2}\left(1-\frac{3}{2}D+\frac{19}{4}D^2-\frac{91}{8}D^3\right)
\left(2x^3+4x^2-6x+5 \right)\nonumber\\
&=&x^3-\frac{5}{2}x^2+\frac{39}{2} x-\frac{169}{4}
\end{eqnarray}

\noindent \textbf{Example 2:} Determine a particular solution of
the equation
\begin{eqnarray}
 y^{\prime\prime\prime}-3 y^{\prime\prime}+2y^{\prime}=x^3-2x^2
\end{eqnarray}
A particular  solution of the above equation  is
\begin{eqnarray}
  Y(x)&=& \frac{1}{D^3-3 D^2+2D}\left(x^3-2x^2 \right)
  =\frac{1}{D(1-D)(2-D)}\left( x^3-2x^2\right)\nonumber\\
  &=&\frac{1}{D(2-D)}\left(1+D+D^2+D^3\right)\left(
  x^3-2x^2\right)\nonumber\\
&=&\frac{1}{D}\left[\frac{1}{2}
\left(1+\frac{1}{2}D+\frac{1}{4}D^2+\frac{1}{8}D^3\right)\right]
(x^3+x^2+2x+2)\nonumber\\
&=&\frac{1}{2}D^{-1}\left(
x^3+\frac{5}{2}x^2+\frac{9}{2}x+\frac{17}{4} \right)
=\frac{1}{8}x^4+\frac{5}{12}x^3+\frac{9}{8}x^2+\frac{17}{8}x
\end{eqnarray}

\subsection{The case $g(x)$ is either a sine or cosine function}

 For $g(x)= A \sin\beta x$ or $g(x)=A\cos\beta x$, the solution is
formally given by Eq.(\ref{formsolu}) as previous cases:
\begin{eqnarray}
Y(x)=\frac{1}{P_n(D)} \left( A \sin\beta x+B \sin\beta x\right)
\label{formtrisolu}
\end{eqnarray}

There are two cases:

\begin{itemize}
\item[(1).] $P_n(D)\sin\beta x\neq 0$ and $P_n(D)\cos\beta x\neq
0$

In this case, we first employ the result of \textbf{Corollary 2}
and the simple algebraic operation $(a+b)(a-b)=a^2-b^2 $ to make
the denominator become a function of $D^2$, and then apply the
eigenvalue substitution rule (\ref{eigensub2}). We use the
following example to illustrate this method.

\textbf{Example:} Determine a particular solution of the equation
\begin{eqnarray}
2y^{\prime\prime\prime}+y^{\prime\prime}-5 y^\prime+3y=3\sin 2x
\label{trigfex1}
\end{eqnarray}
A particular solution of the equation is
\begin{eqnarray}
Y(x)&= &\frac{1}{2D^3+D^2-5D+3}\left(3\sin 2x\right)
=\frac{(2D^3-5D)-(D^2+3)}{(2D^3-5D)^2-(D^2+3)^2}\left(3\sin
2x\right)\nonumber\\
&=&\frac{2D^3-D^2-5D-3}{D^2(2D^2-5)^2-(D^2+3)^2}\left(3\sin
2x\right)\nonumber\\
&=&\frac{3}{(-4)\left[2\cdot(-4)-5\right]^2-(-4+3)^2}\left(
2\frac{d^3}{dx^3}-\frac{d^2}{dx^2}-5\frac{d}{dx}-3\right)\sin
2x\nonumber\\
&=&\frac{3}{677}\left(26\cos 2x-\sin 2x\right)
\end{eqnarray}

\item[(2).] $P_n(D)\sin\beta x= 0$ and $P_n(D)\cos\beta x= 0$

In this case, due to $\displaystyle \left(D^2+\beta^2
\right)\sin\beta x =0$ and $\displaystyle \left(D^2+\beta^2
\right)\cos\beta x =0$, $P_n(D)$ must takes the following form,
\begin{eqnarray}
P_n(D)=P_{n-2k}(D) \left(D^2+\beta^2\right)^k, ~~~ k\geq 1
\end{eqnarray}

 To determine a particular solution in this case, we first
 introduce the following result.

\noindent \textbf{Theorem 7}: There exist
\begin{eqnarray}
 \frac{1}{\left(D^2+\beta^2\right)^{2p}}\cos\beta x
&=& \frac{(-1)^p}{(2p)!(2\beta)^{2p}} x^{2p}\cos\beta x,\nonumber\\
\frac{1}{\left(D^2+\beta^2\right)^{2p}}\sin\beta x
&=& \frac{(-1)^p}{(2p)!(2\beta)^{2p}} x^{2p}\sin\beta x, ~~ p=1,2, \cdots
\label{dsresult1}\\
\frac{1}{\left(D^2+\beta^2\right)^{2p+1}}\cos\beta x
&=& \frac{(-1)^p}{(2p+1)!(2\beta)^{2p+1}} x^{2p+1}\sin\beta x,\nonumber\\
\frac{1}{\left(D^2+\beta^2\right)^{2p+1}}\sin\beta x &=&
\frac{(-1)^{p+1}}{(2p+1)!(2\beta)^{2p+1}} x^{2p+1}\cos\beta x,
~~p=0,1,2, \cdots \label{dsresult2}
\end{eqnarray}

\noindent \textbf{Proof:} Applying the exponential shift rule, we
have
\begin{eqnarray}
&& \left( D^2+\beta^2\right)^k\left( x^k e^{i\beta x}\right)=
e^{i\beta x}\left[ \left(D+i\beta \right)^2+\beta^2\right]^k x^k
\nonumber\\
&&=e^{i\beta x}\left(D+2i\beta\right)^k D^k x^k=e^{i\beta x} i^k
(2\beta)^k k!
\end{eqnarray}

 For $k=2p$, $p=1, 2,\cdots$, using the Euler formula, we obtain
\begin{eqnarray}
\left(D^2+\beta^2 \right)^{2p}\left(x^{2p}\cos\beta x+i
x^{2p}\sin\beta x\right) = (2p)! (2\beta)^{2p} (-1)^p
\left(\cos\beta x+i\sin\beta x\right)
\end{eqnarray}
The real and imaginary parts yield
\begin{eqnarray}
&& \left(D^2+\beta^2 \right)^{2p}\left(x^{2p}\cos\beta x\right)
=(2p)! (2\beta)^{2p} (-1)^p \cos\beta x, \nonumber\\
&& \left(D^2+\beta^2 \right)^{2p}\left(x^{2p}\sin\beta x\right)
=(2p)! (2\beta)^{2p} (-1)^p \sin\beta x
\end{eqnarray}
and lead to the result (\ref{dsresult1}).

Similarly, for $k=2p+1$, $p=0,1,\cdots$, there arises
\begin{eqnarray}
&& \left(D^2+\beta^2 \right)^{2p+1}\left(x^{2p+1}\cos\beta x+i
x^{2p+1}\sin\beta x\right)\nonumber\\
 &=& (2p+1)! (2\beta)^{2p}
(-1)^p \left(i\cos\beta x-\sin\beta x \right)
\end{eqnarray}
The real and imaginary parts give the result (\ref{dsresult2}).

We hence can determine a particular solution (\ref{formtrisolu})
for the case $P_n(D)=0$ as follows:
\begin{eqnarray}
Y(x)&=& \frac{1}{P_{n-2k}(D)
\left(D^2+\beta^2\right)^k}\left(A\cos\beta x+B \sin\beta
x\right)\nonumber\\
&=&
\frac{1}{\left(D^2+\beta^2\right)^k}\left[\frac{1}{P_{n-2k}(D)}\left(A\cos\beta
x+B \sin\beta x\right) \right]
\end{eqnarray}

For example, we determine a particular solution of the
differential equation
\begin{eqnarray}
(D-1)^2(D-2)(D^2+4)^2 y(x)=4 \sin (2x) \label{trigfex2}
\end{eqnarray}
Then
\begin{eqnarray}
Y(x)&=& \frac{1}{(D-1)^2(D-2)(D^2+4)^2} \left(4\sin 2 x\right)
=\frac{4}{(D^2+4)^2}\left[\frac{1}{(D-1)^2(D-2)} \sin
2x\right]\nonumber\\
&=&\frac{4}{(D^2+4)^2}\left[\frac{(D+1)^2
(D+2)}{(D^2-1)^2(D^2-4)}\sin 2x \right]\nonumber\\
&=&-\frac{1}{50}\,\frac{1}{(D^2+4)^2}\left(
\frac{d^3}{dx^3}+4\frac{d^2}{dx^2}+5\frac{d}{dx}+2\right)\sin
2x\nonumber\\
&=&-\frac{1}{50}\,\frac{1}{(D^2+4)^2}\left(2 \cos 2x-14 \sin
2x\right)=\frac{1}{800}x^2\left(\cos 2x-7 \sin 2x\right)
\end{eqnarray}

\end{itemize}

\noindent An alternative method for solving the case $g(x)=A
\cos\beta (x)+B\sin\beta x$ is first using the Euler formula to
express $\cos\beta x$ and $\sin\beta x$  in terms of complex
exponential functions,
\begin{eqnarray}
\cos\beta x =\frac{1}{2}\left(e^{i\beta x}+e^{-i\beta x} \right),
~~\sin\beta x =\frac{1}{2}\left(e^{i\beta x}-e^{-i\beta x} \right)
\end{eqnarray}
and then determining a particular solution in the same way as
$g(x)$ being a complex exponential function. For example, the
particular solutions of the equations (\ref{trigfex1}) and
(\ref{trigfex2}) can be worked out as follows:

\begin{itemize}
\item[(1).] For the equation (\ref{trigfex1}), we have
\begin{eqnarray}
 Y(x)&=& \frac{1}{2D^3+D^2-5D+3}\left(3\sin 2x\right)\nonumber\\
 &=&\frac{1}{2D^3+D^2-5D+3}\left[\frac{3}{2i}\left(e^{2ix}-e^{-2ix}\right)\right]
 \nonumber\\
 &=&\frac{3}{2i}\left[e^{2ix}\left(\frac{1}{2 (D+2i)^3+ (D+2i)^2-5
 (D+2i)+3}\,1\right)\right.\nonumber\\
 && \left.-e^{-2ix}\left(\frac{1}{2 (D-2i)^3+ (D-2i)^2-5 (D-2i)+3}\,1
 \right)\right]\nonumber\\
&=&\frac{3}{2i}\left(\frac{e^{2ix}}{-1-26 i}-\frac{e^{-2ix}}{-1+26
i} \right)\nonumber\\
&=& \frac{3}{2i}\frac{1}{677}\left[(-1+26i)(\cos 2x+i \sin
2x)-(-1-26i)(\cos 2x-i \sin 2x)\right]\nonumber\\
&=&\frac{3}{677}\left( 26 \cos 2x-\sin 2x\right)
\end{eqnarray}

\item[(2).] For the equation (\ref{trigfex2}), we have
\begin{eqnarray}
 Y(x)&=& \frac{1}{(D-1)^2(D-2)(D^2+4)^2}\left(4\sin 2x\right)\nonumber\\
 &=&\frac{1}{(D-1)^2(D-2)(D^2+4)^2} \left[\frac{4}{2i}
 \left(e^{2ix}-e^{-2ix}\right)\right]
 \nonumber\\
 &=&\frac{2}{i}\left[e^{2ix}
 \left(\frac{1}{ (D+2i-1)^2 (D+2i-2)[(D+2i)^2+4]^2}\,1\right)\right.\nonumber\\
 && \left.-e^{-2ix}\left(\frac{1}{ (D-2i-1)^2 (D-2i-2)[(D-2i)^2+4]^2}\,1\right)
 \right]\nonumber\\
&=&\frac{2}{i}\left[\frac{e^{2ix}}{(2i-1)^2(2i-2)}\left(\frac{1}{D^2(D+4i)^2}\,1\right)
-\frac{e^{-2ix}}{(-2i-1)^2(-2i-2)}\left(\frac{1}{D^2(D-4i)^2}\,1\right)
\right]\nonumber\\
&=&\frac{1}{16i}\left[  \frac{e^{2ix}}{(1-2i)^2 (1-i)}+
\frac{e^{-2ix}}{(1+2i)^2 (1+i)} \right] D^{-2}\,1\nonumber\\
&=&\frac{1}{800i}\left[(-3+4i)(1+i)(\cos 2x+i\sin 2x)
+(-3-4i)(1-i)(\cos 2x-i\sin 2x)\right]\frac{1}{2}x^2
\nonumber\\
 &=&\frac{1}{800}x^2\left(\cos 2x-7\sin 2x  \right)
\end{eqnarray}

\end{itemize}

\subsection{The case $g(x)$ is a product of exponential
function and polynomial function}

\noindent For $g(x)=P_m(x) e^{\alpha x}$, one should first apply
the exponential shift rule to move out of the exponential
function, and then proceed to determine a particular solution as
the case $g(x)$ being a polynomial function:
\begin{eqnarray}
Y(x)=\frac{1}{P_n(D)}\left( P_m(x)e^{\alpha x}\right) = e^{\alpha
x}\left[\frac{1}{P_n(D+\alpha)} P_m(x)\right]
\end{eqnarray}
For examples,

\begin{itemize}
\item[(1).]  A particular solution of the equation
\begin{eqnarray}
(D-3)^2 (D^2-2D+5)(D+2) y(x)= (x^2-3x+1)e^{2x}
\end{eqnarray}
is
\begin{eqnarray}
 Y(x)&=& \frac{1}{(D-3)^2 [(D-1)^2+4] (D+2)}
 \left[(x^2-3x+1)e^{2x}\right]\nonumber\\
 &=& e^{2x}\frac{1}{(D-1)^2 [(D+1)^2+4] (D+4)}(x^2-3x+1)\nonumber\\
&=&e^{2x}\frac{1}{20-27D+2 D^3+4D^4+D^5}(x^2-3x+1)\nonumber\\
&=& e^{2x}\frac{1}{20-27D}(x^2-3x+1)=e^{2x}\frac{1}{20} \left[
1+\frac{27}{20}\frac{d}{dx}+\left(\frac{27}{20}\right)^2\frac{d^2}{dx^2}\right]
(x^2-3x+1)\nonumber\\
&=&
\frac{1}{20}\left(x^2-\frac{3}{10}x+\frac{119}{200}\right)e^{2x}
\end{eqnarray}

\item[(2).]  A particular solution of the equation
\begin{eqnarray}
(D-3) (D-2)^2(D+1) y(x)= (4x-2)e^{2x}
\end{eqnarray}
is
\begin{eqnarray}
 Y(x)&=& \frac{1}{(D-2)^2(D-3)(D+1)}
 \left[(4x-2)e^{2x}\right]= e^{2x}\frac{1}{(D-1)(D+3)D^2}(4x-2)\nonumber\\
&=&e^{2x}\frac{1}{D^2(D^2+2D-3)}(4x-2)
=e^{2x}\frac{1}{D^2}\left[\frac{1}{2D-3}(4x-2)\right]\nonumber\\
&=&
e^{2x}\frac{1}{D^2}\left[-\frac{1}{3}
\left(1+\frac{2}{3}\frac{d}{dx}\right)(4x-2)\right]
=-\frac{1}{3}e^{2x}D^{-2}\left(4x+\frac{2}{3}\right)\nonumber\\
&=& -\frac{1}{9}x^2 (2x+1)e^{2x}
\end{eqnarray}

\end{itemize}

\subsection{The case $g(x)$ is a product of polynomial function
and sine or cosine functions}

\noindent In this case, $g(x)=P_m(x) (A\cos \beta x+B \sin\beta
x)$. The most convenient way of determining a particular solution
is to first express sine and cosine function in terms of complex
exponential functions, and then use the method of dealing with the
product of a polynomial function and exponential function to find
a particular solution:
\begin{eqnarray}
&& Y(x)= \frac{1}{P_n(D)}\left[ P_m(x) \left(A\cos \beta x+B
\sin\beta x \right)\right] \nonumber\\
&=& \frac{1}{2}\left(A-iB\right) e^{i\beta x}
\left[\frac{1}{P_n(D+i\beta)}
P_m(x)\right]+\frac{1}{2}\left(A+iB\right) e^{-i\beta x}
\left[\frac{1}{P_n(D-i\beta)} P_m(x)\right]
\end{eqnarray}

\noindent For examples,
\begin{itemize}
\item[(1).] A particular solution to the differential equation
\begin{eqnarray}
y^{\prime\prime}-4y=  (x^2-3)\sin 2x
\end{eqnarray}
is
\begin{eqnarray}
Y(x)&=&\frac{1}{D^2-4}\left[\left(x^2-3\right)\sin 2x  \right] =
\frac{1}{D^2-4}\left[\left(x^2-3\right)\frac{1}{2i}\left(e^{2ix}-e^{-2ix}\right)
\right]\left(x^2-3\right)\nonumber\\
&=&\frac{1}{2i}\left[
e^{2ix}\frac{1}{(D+2i)^2-4}-e^{-2ix}\frac{1}{(D-2i)^2-4}\right]\left(x^2-3\right)
\nonumber\\
&=&\frac{1}{2i}\left[
e^{2ix}\frac{1}{D^2+4iD-8}-e^{-2ix}\frac{1}{D^2-4iD-8}\right]\left(x^2-3\right)
\nonumber\\
&=&\frac{i}{16}\left[e^{2ix}\left(1+\frac{1}{2}iD-\frac{1}{8}D^2\right)
-e^{-2ix}\left(1-\frac{1}{2}iD-\frac{1}{8}D^2\right) \right]
(x^2-3)\nonumber\\
&=& -\frac{1}{32}\left[ \left(4x^2-13\right)\sin 2x +4x \cos
2x\right]
\end{eqnarray}


\item[(2).] A particular solution to the equation
\begin{eqnarray}
y^{\prime\prime}+4y= 4x^2 \cos 2x
\end{eqnarray}
is
\begin{eqnarray}
Y(x)&=&\frac{1}{D^2+4}\left(4x^2\cos 2x
\right)=\frac{1}{D^2+4}\left[2x^2\left(e^{2ix}+e^{-2ix}\right)
\right]\nonumber\\
&=&2 \left[e^{2ix}\frac{1}{(D+2i)^2+4}x^2+
e^{-2ix}\frac{1}{(D-2i)^2+4}x^2\right]\nonumber\\
&=&2\left[e^{2ix}\frac{1}{D(D+4i)}x^2+
e^{-2ix}\frac{1}{D(D-4i)}x^2\right]\nonumber\\
&=&\frac{1}{4i}\left[e^{2ix}D^{-1}\left(2
x^2+ix-\frac{1}{4}\right) -e^{-2ix}D^{-1}\left(2
x^2-ix-\frac{1}{4}\right)\right]\nonumber\\
&=&\frac{1}{24}\left[6x^2\cos 2x+x (8x^2-3)\sin 2x \right]
\end{eqnarray}


\end{itemize}

\subsection{The case $g(x)$ is a product of exponential function
and sine or cosine functions}

\noindent For $g(x)=e^{\alpha x} (A\cos \beta x+B \sin\beta x)$,
as other cases involving sine or cosine functions, the best way of
determining a
 particular solution for an equation is to first express sine
and cosine function in terms of complex exponential functions.
There are two possibilities.

\begin{itemize}
\item[(1).] If $P_n(\alpha\pm i\beta)\neq 0$, then one can
directly apply the eigenvalue substitution rule to get a
particular solution:
\begin{eqnarray}
Y(x)&=&\frac{1}{P_n(D)}\left[e^{\alpha x}\left(A\cos \beta x+B
\sin\beta x\right)
\right]\nonumber\\
&=& \frac{A-iB}{2P_n(\alpha+i\beta)}e^{(\alpha+i\beta)x} +
\frac{A+iB}{2P_n(\alpha-i\beta)}e^{(\alpha-i\beta)x}
\end{eqnarray}

For example, a particular solution of the equation
\begin{eqnarray}
y^{\prime\prime}-2y^{\prime}+2y=e^{2x} \left(2\cos x-6\sin
x\right)
\end{eqnarray}
is
\begin{eqnarray}
Y(x)&=&\frac{1}{D^2-2D+2}\left[ e^{2x} \left(2\cos x-6\sin
x\right)\right] \nonumber\\
&=&\frac{1}{D^2-2D+2}\left[(1+3
i)e^{(2+i)x}+(1-3i)e^{(2-i)x} \right]\nonumber\\
&=&(1+3i)\frac{e^{(2+i)x}}{(2+i)^2-2(2+i)+2}
+(1-3i)\frac{e^{(2-i)x}}{(2-i)^2-2(2-i)+2}\nonumber\\
&=&\frac{2}{5}e^{2x}\left(7\cos x-\sin x \right)
\end{eqnarray}

An alternative way is to first utilize the exponential rule to
move out $e^{\alpha x}$ and then use the method of dealing with
sine or cosine function to determine a particular solution. For
example,  a particular solution of the above equation in this
approach is obtained as follows:
\begin{eqnarray}
Y(x)&=&e^{2x}\frac{1}{(D+2)^2-2(D+2)+2} \left(2\cos x-6\sin
x\right) = e^{2x}\frac{1}{D^2+2D+2}\left(2\cos x-6\sin
x \right)\nonumber\\
&=&e^{2x}\frac{D^2+2-2D}{(D^2+2)^2-4D^2}\left(2\cos x-6\sin x
\right) =\frac{1}{5}e^{2x}\left(\frac{d^2}{dx^2}-2
\frac{d}{dx}+2\right)\left(2\cos x-6\sin
x \right)\nonumber\\
&=&\frac{2}{5}e^{2x}\left(7\cos x-\sin x \right)
\end{eqnarray}

\item[(2).] If $P_{n}(\alpha\pm i\beta)=0$, one should first apply
the exponential shift rule and then the definition for the inverse
of differential operator. For example, we determine a  particular
solution of the equation
\begin{eqnarray}
(D-1)(D^2-2D+5)y=4e^{3x}\cos 2x
\end{eqnarray}
Then
\begin{eqnarray}
Y(x)&=&\frac{1}{(D-1)(D^2-6D+13)}\left(4e^{3x}\cos
2x\right)\nonumber\\
&=&\frac{1}{[(D-3)^2+4](D-1)}\left\{2\left[e^{(3+2i)x}+e^{(3-2i)x}\right\}
\right]\nonumber\\
&=&\frac{1}{(D-3)^2+4}\left[ \frac{2}{2+2i}e^{(3+2i)x}
 \right]+\frac{1}{(D-3)^2+4}\left[ \frac{2}{2-2i}e^{(3-2i)x}
 \right]\nonumber\\
 &=&\frac{1}{2}\left(1-i\right)e^{(3+2i)x}\frac{1}{D(D+4i)}\,1
 +\frac{1}{2}\left(1+i\right)e^{(3-2i)x}\frac{1}{D(D-4i)}\,1
 \nonumber\\
 &=&\frac{1}{8i}xe^{3x}(1-i)(\cos 2x+i\sin 2x)
 -\frac{1}{8i}xe^{3x}(1+i)(\cos 2x+i\sin 2x)\nonumber\\
 &=&\frac{1}{4}x e^{3x}\left(\sin 2x -\cos 2x \right)
\end{eqnarray}

One can also work first take out $e^{\alpha x}$ with the
exponential shift rule and then apply the results
(\ref{dsresult1}) and (\ref{dsresult2}) of \textbf{Theorem 7}. The
particular solution of the above example in this approach reads
\begin{eqnarray}
Y(x)&=&\frac{1}{(D-1)[(D-3)^2+4])}\left(4e^{3x}\cos 2x\right)
=4e^{3x}\frac{1}{(D+2)(D^2+4)}\cos 2x \nonumber\\
&=&4e^{3x}\frac{1}{D^2+4}\left(\frac{D-2}{D^2-4}\cos 2x
\right)=e^{3x}\frac{1}{D^2+4}\left(\sin 2x+\cos 2x
\right)\nonumber\\
&=&\frac{1}{4}xe^{3x}\left(\sin 2x-\cos 2x \right)
\end{eqnarray}

\end{itemize}

\subsection{The case $g(x)$ is a product of polynomial function, exponential function
and sine or cosine functions}

\noindent Finally, we consider the most general case, $g(x)=P_m(x)
e^{\alpha x} (A\cos\beta x+B\sin\beta x)$. The method of getting a
particular is the same as the previous case that $g(x)$ is a
product of a polynomial function and an exponential function:
\begin{eqnarray}
Y(x)&=&\frac{1}{P_n(D)}\left[P_m(x)e^{\alpha x}\left(A \cos\beta
x+B\sin \beta x \right) \right]\nonumber\\
&=&\frac{1}{2}\left(A-iB\right)\frac{1}{P_n(D+\alpha+i\beta)}P_m(x)
+\frac{1}{2}\left(A+iB\right)\frac{1}{P_n(D+\alpha-i\beta)}P_m(x)
\end{eqnarray}

\noindent For examples,
\begin{itemize}
\item[(1).] A particular solution to the differential equation
\begin{eqnarray}
y^{\prime\prime}-5y^{\prime}+6y=e^x\cos 2x+e^{2x}\left(3x+4
\right)\sin x
\end{eqnarray}
is
\begin{eqnarray}
Y(x)&=&\frac{1}{D^2-5D+6}\left[e^x\cos 2x+e^{2x}\left(3x+4
\right)\sin x\right]=Y_1(x)+Y_2(x),\\
Y_1(x)&=&\frac{1}{D^2-5D+6}\left(e^x\cos 2x
\right)=\frac{1}{(D-2)(D-3)}\left(e^x\cos 2x \right)\nonumber\\
&=& e^x\frac{1}{(D-1)(D-2)}\cos
2x=e^x\frac{(D+1)(D+2)}{(D^2-1)(D^2-4)} \cos 2x\nonumber\\
&=& \frac{1}{40}e^x \left(\frac{d}{dx}+1
\right)\left(\frac{d}{dx}+2 \right)\cos 2x=-\frac{1}{20}e^x \left(
\cos 2x+3 \sin 2x\right),\\
Y_2(x)&=&\frac{1}{D^2-5D+6}\left[\left(3x+4\right)e^{2x}\cos x
\right]\nonumber\\
&=&
\frac{1}{2}\,\frac{1}{(D-2)(D-3)}\left[(3x+4)\left(e^{(2+i)x}+e^{(2-i)
x}\right) \right]\nonumber\\
&=&\frac{1}{2}e^{(2+i)x}\frac{1}{(D+i)(D-1+i)}\left(3x+4\right)\nonumber\\
&&
+\frac{1}{2}e^{(2-i)x}\frac{1}{(D-i)(D-1-i)}\left(3x+4\right)\nonumber\\
&=&\frac{1}{2i}e^{(2+i)x}\frac{1}{D-1+i}\left( 3x+4+3i\right)
 -  \frac{1}{2i}e^{(2-i)x}\frac{1}{D-1-i}\left( 3x+4-3i\right) \nonumber\\
&=&-\frac{1}{2}e^{2x}\left[(3x+10) \cos x+(3x+1)\sin x \right]
\end{eqnarray}

\item[(2).] A particular solution of the equation
\begin{eqnarray}
y^{\prime\prime}+2y^{\prime}+2y=e^{-x}\left(3+2\sin x+4 x^2\cos x
\right)
\end{eqnarray}
is
\begin{eqnarray}
Y(x)=\frac{1}{D^2+2D+2}\left[e^{-x}\left(3+2\sin x+4 x^2\cos x
\right)\right]=Y_1(x)+Y_2(x)+Y_3(x)
\end{eqnarray}
where
\begin{eqnarray}
Y_1(x)&=&\frac{1}{D^2+2D+2}\left(3e^{-x}\right)=\frac{3e^{-x}}{(-1)^2-2+2}=3e^{-x},
\nonumber\\
Y_2(x) &=&\frac{1}{D^2+2D+2}\left(2e^{-x}\sin x\right) =
\frac{1}{(D+1)^2+1}\left( 2e^{-x}\sin
x\right)\nonumber\\
&=& 2e^{-x}\frac{1}{D^2+1}\sin x=-xe^{-x}\cos x,\\
Y_3(x) &=&\frac{1}{D^2+2D+2}\left(4x^2e^{-x}\cos x\right) =
\frac{1}{(D+1)^2+1}\left( 4x^2e^{-x}\cos x\right)\nonumber\\
&=&\frac{1}{(D+1)^2+1}\left[2 x^2\left(e^{(-1+i)x}+e^{(-1-i)x}
\right)\right]\nonumber\\
&=&2e^{(-1+i)x}\frac{1}{D(D+2i)}x^2+2e^{(-1-i)x}\frac{1}{D(D-2i)}x^2\nonumber\\
&=&-i e^{(-1+i)x} D^{-1}\left(x^2+i x-\frac{1}{2} \right)
+ie^{(-1+i)x} D^{-1}\left(x^2-i x-\frac{1}{2} \right)\nonumber\\
&=&\frac{1}{3}e^{-x}\left[ 3x^2\cos x+x (2x^2-3)\sin x\right]
\end{eqnarray}

\end{itemize}

\section{Summary}

\noindent We have introduced both the mathematical principle and
computing technique of using differential operator method  to find
a particular solution of an ordinary nonhomogeneous linear
differential equation with constant coefficients when the
nonhomogeneous term is a polynomial function, exponential
function, sine function, cosine function or any possible product
of these functions. We have reviewed some rules in the
differential operations including the eigenvalue substitution rule
and the exponential shift rule. Furthermore, we have used a number
of examples to illustrate the practical application of this
approach. In particular, comparing with the introductions on this
differential operator method scattered in some lecture notes, we
propose and highlight the application of the inverse of
differential operator in determining a particular solution,  and
overcome the difficulty when a polynomial of differential operator
in a differential equation is singular. The differential operator
method has great advantages over the well-known method of
undetermined coefficients introduced in textbooks in determining a
particular solution. Therefore, it is meaningful  to make a
systematic introduction and put an emphasis on this method. We
hope that this method can be systematically introduced in
textbooks and widely used for determining a particular solution of
an ordinary  nonhomogeneous linear differential equation with
constant coefficients, in parallel to the method of undetermined
coefficients.

\acknowledgments \noindent I would like to thank my Colleagues in
the Group of Applied Mathematics at SUNY Polytechnic Institute for
useful remarks.  I am especially indebted to Dr. Andrea Dziubek
and Dr. Edmond Rusjan for various discussions and encouragements.

\end{document}